\documentclass[12pt]{amsart}

\usepackage{amssymb,pstricks}
\usepackage[all]{xy}
\begin{document}
\newtheorem{cor}{Corollary}
\newtheorem{theorem}[cor]{Theorem}
\newtheorem{prop}[cor]{Proposition}
\newtheorem{lemma}[cor]{Lemma}
\theoremstyle{definition}
\newtheorem{defi}[cor]{Definition}
\theoremstyle{remark}
\newtheorem{remark}[cor]{Remark}
\newtheorem{example}[cor]{Example}

\newcommand{\dem}{{\noindent{\sc {Proof. }}}}
\newcommand{\bX}{\partial X}
\newcommand{\cC}{\mathcal{C}}
\newcommand{\cD}{\mathcal{D}}
\newcommand{\cE}{\mathcal{E}}
\newcommand{\cF}{\mathcal{F}}
\newcommand{\cG}{\mathcal{G}}
\newcommand{\cH}{\mathcal{H}}
\newcommand{\cI}{\mathcal{I}}
\newcommand{\cN}{\mathcal{N}}
\newcommand{\cV}{\mathcal{V}}
\newcommand{\ccV}{{}^c\cV}
\newcommand{\coker}{\operatorname{coker}}
\newcommand{\cssX}{{}^c\!S^*\oX}
\newcommand{\ctX}{{}^c\!T\oX}
\newcommand{\ctsX}{{}^c\!T^*\oX}
\newcommand{\cun}{\cC^{\infty}}
\newcommand{\cz}{{\mathbb C}}
\newcommand{\dc}{\cD_c}
\newcommand{\en}{{\epsilon_n}}
\newcommand{\Hom}{\operatorname{Hom}}
\newcommand{\htX}{{}^h\!TX}
\newcommand{\ind}{\operatorname{index}}
\newcommand{\iX}{{\stackrel{\circ}{X}}}
\newcommand{\nz}{{\mathbb N}}
\newcommand{\oX}{\overline{X}}
\newcommand{\pc}{\Psi_c}
\newcommand{\ps}{\Psi_\sus}
\newcommand{\px}{\partial_x}
\newcommand{\re}{\operatorname{Re}}
\newcommand{\Res}{\operatorname{Res}}
\newcommand{\rc}{{\rho_c}}
\newcommand{\rz}{{\mathbb R}}
\newcommand{\scal}{\operatorname{scal}}
\newcommand{\sD}{\sigma(D)}
\newcommand{\sign}{\operatorname{sign}}
\newcommand{\Spec}{\operatorname{Spec}}
\newcommand{\Spin}{{\operatorname{Spin}}}
\newcommand{\sssM}{{S_\sus^*(M)}}
\newcommand{\sus}{{\operatorname{sus}}}
\newcommand{\ta}{{\tilde{a}}}
\newcommand{\tD}{{\tilde{D}}}
\newcommand{\td}{{\tilde{\delta}}}
\newcommand{\tg}{\tilde{g}}
\newcommand{\tr}{\operatorname{tr}}
\newcommand{\Tr}{\operatorname{Tr}}
\newcommand{\tssM}{T_\sus^* M}
\newcommand{\tX}{{\tilde{X}}}
\newcommand{\tY}{{\tilde{Y}}}
\newcommand{\Vol}{\operatorname{Vol}}

\title{Weyl laws on open manifolds}
\author{Sergiu Moroianu}
\thanks{The author has been partially supported by the
Research and Training Network HPRN-CT-1999-00118 ``Geometric Analysis'' 
funded by the European Commission}
\date{\today}
\address{Institutul de Matematic\u{a} al Academiei Rom\^{a}ne\\
P.O. Box 1-764\\RO-70700
Bucharest, Romania}
\address{Universit\'e Paul Sabatier, UFR MIG,
118 route de Narbonne, 31062 Toulouse, France}
\email{moroianu@alum.mit.edu}
\begin{abstract}
Under suitable invertibility hypothesis, the spectrum 
of the Dirac operator on certain open spin Riemannian 
manifolds is discrete, and 
obeys a growth law depending qualitatively on the 
(in)finiteness of the volume.
\end{abstract}
\maketitle

\section{Introduction} 

Let $M$ be a closed Riemannian manifold, 
$\cE\to M$ a Hermitian vector bundle 
and $D:\cun(M,\cE)\to\cun(M,\cE)$ an 
elliptic, symmetric, positive differential operator of order $k>0$. 
Then the spectrum of $D$ is 
discrete and the eigenvalues accumulate towards infinity obeying the 
Weyl law:
\begin{equation}\label{weyl}
\lim_{\lambda\to\infty}\frac{N(\lambda)}{\lambda^{\dim(M)/k}}= C
\end{equation}
where $N(\lambda):=\#\{s\in\Spec(D); |s|<\lambda\}$ 
is the counting function and $C$ is an explicit constant depending on $\dim(M)$ 
and on the principal symbol of $D$.
This fact is classically proved using the heat trace expansion \cite{mipl} 
when $D$ is a Laplacian. Note that one can obtain 
better estimates of the remainder in (\ref{weyl}), see \cite{hor}.

In this paper we derive similar laws for Dirac 
operators on certain open spin Riemannian manifolds. 

If the manifold $M$ is not compact, not much can be said in general
about the spectrum of $D$. However, if $M$ is complete then the Laplacian on 
forms and the Dirac operator are essentially self-adjoint, and their essential
spectra may be non-empty. Several results are also known for compact, 
incomplete manifolds with boundary. The original result of Weyl 
states that (\ref{weyl}) holds
for the Dirichlet Laplacian
on a compact domain with smooth boundary in $\rz^n$. 
Let now $\oX$ be a compact manifold with boundary $M$, and $X$ its interior.
Then (\ref{weyl}) holds \cite{aps1} for the Dirac operator on $X$ endowed 
with a product-type metric near $M$ and a suitable non-local boundary condition.
The results of \cite{aps1} follow from an 
explicit computation on the cylinder, which breaks down for warped product
metrics. A significant progress in this context was brought by Melrose
\cite{melaps}, who computed the index of the Atiyah-Patodi-Singer 
non-local boundary value problem as an $L^2$ index, using his b-calculus on 
the complete manifold obtained from $\oX$ by gluing infinite cylinders.
Unlike the APS operator, the associated b-operator has non-empty essential 
spectrum.

The present work is motivated by a result of B\"ar
\cite{baer}, who showed that the essential spectrum of the Dirac operator on
complete spin hyperbolic manifolds of finite volume is either empty or the 
whole real line. Thus, let $\oX$ be a smooth $n$-dimensional compact manifold 
with closed boundary $M$, and $x$ a boundary-defining function. This means 
\begin{itemize}
\item $x\in\cun(\oX,[0,\infty))$;
\item $M=\{x=0\}$;
\item $dx$ never vanishes at $x=0$.
\end{itemize}
Let $g_0$ be a \emph{cusp metric} on $\oX$, i.e., a 
(complete) Riemannian metric on $X:=\oX\setminus M$ 
which in local coordinates near the boundary takes the form
\begin{equation}\label{cume}
g_0=a_{00}(x,y)\frac{dx^2}{x^4}+\sum_{j=1}^{n-1} a_{0j}(x,y)
\frac{dx}{x^2}dy_j+\sum_{i,j=1}^{n-1}a_{ij}(x,y)dy_idy_j
\end{equation}
such that the matrix $A=(a_{\alpha\beta})$ is smooth and non-degenerate
down to $x=0$. Let $h:=\sum_{i,j=1}^{n-1}a_{ij}(0,y)dy_idy_j$  
denote the induced Riemannian metric on $M$.
For simplicity, in this introduction we present our results 
under the additional assumption  
\begin{equation}\label{adias}
{a_{00}}_{|M}\equiv 1,\ \ {a_{0j}}_{|M}\equiv 0.
\end{equation}
Following Melrose \cite{melaps}, we call $g_0$ satisfying (\ref{adias})
an \emph{exact cusp metric}.
We are interested in the Riemannian metrics
\begin{equation}\label{gp}
g_p:=x^{2p}g_0
\end{equation}
for $p\in\rz$. Important particular cases are obtained when
$g_0$ is a cylindrical metric on $(0,\epsilon)\times M\subset X$, $g_0=x^{-4}dx^2+h$. 
By the change of variable $y^q=x$, $q>0$
(which changes the smooth structure of $\oX$
but not that of $X$) and for varying $p\in\rz$ we can get any metric
of the type
$$x^adx^2+x^b h$$
with $a,b\in\rz$ provided $b-a> 2$, in particular $g_p$ cannot 
be conformally conical. Such metrics include all metric horns \cite{lepe}, 
and complete hyperbolic manifolds of finite volume.
The above change of variable is rather artificial since it destroys 
smoothness of coefficients in the general case (\ref{cume}). 

We first study the self-adjointness of the Dirac operator of the metric $g_p$
acting on spinors.

\begin{theorem}\label{th1}
Let $X$ have a spin structure such that the Dirac operator 
$D^h$ for the induced spin structure on $(M,h)$ is invertible.
Assume that $g_0$ satisfies (\ref{adias}).
Then the Dirac operator $D_p$ with domain $\cun_c(X,\Sigma)$
is essentially self-adjoint in $L^2(X,\Sigma,g_p)$.
\end{theorem}

This result is standard for $p\leq 1$ even without the invertibility 
hypothesis, since then $g_p$ is complete. For $p>1$
and $g_0=x^{-4}dx^2+h$ we recover a result from \cite{lepe} through the
change of variable $y=x^{p-1}$.
For us the most interesting phenomena will occur for $0<p\leq 1/n$.

It is important to understand how restrictive the invertibility hypothesis 
really is. On one hand, it is conjectured that any spin manifold $M$ admits 
metrics with harmonic spinors, provided $\dim(M)\geq 3$. On the other hand, 
it is also conjectured that (for $M$ connected)
for \emph{generic} metrics the dimension
of the space of harmonic spinors is the minimal dimension prescribed by 
index theory. Parts of this conjecture were proved by Maier \cite{maier}
and B\"ar and Dahl \cite{badh}. In our case, by the cobordism 
invariance of the index, the minimal dimension is $0$.
Thus, for $M$ connected, the invertibility hypothesis is automatically 
satisfied for generic metrics if 
$\dim(M)\leq4$, or if $\pi_1(M)$ vanishes or is cyclic of odd order.
An indirect way of ensuring invertibility 
of $D^h$ is (following Lichnerowicz) asking that $\scal_h\geq 0$, with 
strict inequality at least at one point.

Incidentally, under the hypothesis of Theorem \ref{th1}, 
$D_p^+$ is in fact Fredholm, and for $\dim X$ even,  
$$\ind(D_p^+)=\int_{\oX} \hat{A}(g_p)-\eta(D^h)/2.$$
The $\hat{A}$ form is a conformal invariant, hence, like $\ind(D_p^+)$,
it is independent of $p$. We will not pursue this fact further.

Next, we look for conditions on $p$ to ensure the vanishing of 
the essential spectrum of $D_p$. 

\begin{theorem}\label{th2}
Under the hypothesis of Theorem \ref{th1}, the operator $D_p$
has pure point spectrum if and only if $p>0$.
\end{theorem}
This greatly generalizes the above cited result of B\"ar \cite{baer}.

%

Our main result describes the asymptotics of the counting 
function for the eigenvalues of $D_p$.

\begin{theorem}\label{th3}
Let $\oX$ be a compact $n$-dimensional 
spin manifold with boundary $M$, $g_p$ the metric on 
$X$ given by (\ref{gp}) and $D_p$ the associated Dirac operator on spinors. 
Assume that  (\ref{adias}) holds and that the spin structure on $\oX$
induces an invertible Dirac operator on $(M,h)$. 
For $p>0$ let $N_p(\lambda)$ denote the counting function
for the eigenvalues of $D_p$ (well-defined by Theorem \ref{th2}). Then
\begin{itemize}
\item For $1/n< p<\infty$, 
$$N_p(\lambda) \approx \lambda^n 
\tfrac{\Vol(X,g_p)\Vol(S^{n-1})2^{[n/2]}}{n(2\pi)^n}.$$
\item For $p=1/n$,
$$N_p(\lambda) \approx\lambda^n\log \lambda
\tfrac{\Vol(M,h)\Vol(S^{n-1})2^{[n/2]}}{(2\pi)^n}.$$
\item For $0<p<1/n$, 
$$N_p(\lambda) \approx \lambda^{1/p}
\tfrac{\Gamma\left(\tfrac{1-p}{2p}\right)\zeta\left(D^h, 
\tfrac{1}{p}-1\right)\en}
{2\sqrt{\pi}\Gamma\left(\tfrac{1}{2p}\right)},$$
where $\en$ equals $1$ if $n$ is odd, $2$ if $n$ is even.
\end{itemize}
\end{theorem}
Here $\Gamma$ is the Gamma function and $\zeta(D^h,z):=
\Tr\left((D^h)^2\right)^{-z/2}$
is the zeta function of the absolute value of $D^h$, which is 
holomorphic for $\Re(z)>n-1=\dim(M)$.

It is interesting to note that $\Vol(X,g_p))$ is finite exactly when 
it arises as coefficient in the asymptotic 
law. In that case the formula bears no difference from the case 
of closed manifolds. The critical conformal factor 
$p=1/n$ marks the start of a remarkable change in the growth rate. 
As $p$ approaches $0$, the metric becomes close to the 
asymptotically cylindrical metric, for which the essential spectrum 
does not vanish. This is reflected by the growth rate 
becoming infinite. Remarkably, some sort of spectral asymptotics
for the Laplacian survive for $g_0$, as shown by Christiansen 
and Zworski \cite{christiansen}.

Weyl asymptotics may be deduced by two related approaches -- 
heat trace expansions, via the beautiful and elementary
result of Karamata, and from meromorphic extensions of zeta functions via
the Ikehara theorem.
Complex powers are weaker but more conceptual than the
heat kernel, in particular they can be constructed in wide generality
\cite{alnv} following the method of Guillemin \cite{gui} and Bucicovschi 
\cite{buc}.
Thus our strategy of proof will be to construct the complex powers of the square of 
our Dirac operators inside a calculus of pseudodifferential operators, and 
examine the leading pole of the associated zeta function. Rather than 
constructing a different calculus for each $p$, which would be highly 
problematic, we rely on Melrose's \emph{cusp calculus}, 
corresponding to $p=0$ (note that geometrically it does not correspond 
to cuspidal but to cylindrical ends). For each $p$ we 
find an isometry between $L^2(\Sigma_p)$ and $L^2(\Sigma_0)$ under which $D_p$
transforms to a cusp operator of order $(1,p)$. The zeta function of 
this operator is a meromorphic function, with first pole at
$\max\{n,1/p\}$, double if $n=1/p$ and simple otherwise. We use then 
the Ikehara theorem \cite{ikehara} 
and its generalization by Delange \cite{delange}
to deduce Theorem \ref{th3}. The invertibility of $D^h$
means that the operator $D_p$ is fully elliptic, hence determines
the domain of its closure and yields Theorem \ref{th1}. As for Theorem \ref{th2},
the condition $p>0$ is equivalent to the parametrix of $D_p$ 
inside the cusp calculus being compact.

It might be possible to retrieve our results by exhibiting small-time
asymptotics of the heat kernel of $D_p^2$. Paul Loya (see e.g., \cite{loya}) has 
constructed such heat kernels and studied their asymptotics in a 
variety of non-standard situations. Such asymptotics would be stronger than the 
meromorphy of the zeta function; in exchange, they would 
yield immediately Theorem \ref{th3} via Karamata's theorem.

Let us mention a third possible approach, i.e., the asymptotics of the wave
trace. Even the first step of identifying this trace with a tempered 
distribution depends mildly on some \emph{a priori} knowledge of 
Weyl asymptotics. Nevertheless, this approach would likely give also 
an estimate of the error term in Theorem \ref{th3}.

We finally mention without proof that
our results seem to extend to manifolds with corners of
codimension $d\geq 2$. For various conformal factors,
one gets as possible growth rates $\lambda^n\log^k \lambda$, 
$k=0,\ldots, d$, and also
$\lambda^{1/p}\log^k \lambda$, $0<p<1/n$, $k=0,\ldots,d-1$.

{\small
{\bf Acknowledgments.}
I wish to thank Andrei Moroianu for several useful discussions.}

\section{The cusp structure}\label{section2}

Consider a general cusp metric given by (\ref{cume}) without the condition
(\ref{adias}), and the closely related metrics on $X$ obtained by the 
conformal transformation
of weight $x^p, p\in\rz$. We call the metric $g_p=x^{2p}g_0$ a 
\emph{p-metric}. Note that the volume of $g_p$ is finite if and only if $p>1/n$.

\begin{example}
A complete hyperbolic manifold $X$ of finite volume is isometric (outside
a compact set) to a warped product cylinder
$$((-\infty,0]\times M,dt^2+e^{-2t}h)$$
where $M$ may be disconnected. Moreover, $h$ is flat and independent of $t$.
By the change of variables $x=e^{-t}$ near $-\infty$, the hyperbolic 
metric transforms to $x^2h+\frac{dx^2}{x^2}$ near $x=0$. Thus a 
hyperbolic metric of finite volume is a $1$-metric on the compactification
$\oX:=X\cup\{-\infty\}\times M$. 
\end{example}

Let 
$$\ccV(\oX):=\{V\in\cun(\oX,T\oX); dx(V)\in x^2\cun(\oX,T\oX)\}\subset\cV(\oX)$$
be the Lie sub-algebra of \emph{cusp vector fields} on $\oX$. 
In local coordinates $(x,y_j)$
near $x=0$, a vector field in $\ccV(\oX)$ can be uniquely written
as a linear combination 
$$a(x,y)x^2\partial_x+\sum_{j=1}^{n-1}b_j(x,y)\partial_{y_j}$$
with smooth coefficients $a,b_j$.
Thus $\ccV(\oX)$ is a locally free $\cun(\oX)$-module. 
By the Serre-Swan
theorem there exists a vector bundle $\ctX\to \oX$ such that
$\ccV(\oX)=\cun(\oX,\ctX)$. Moreover, the inclusion $\ccV(\oX)\hookrightarrow
\cun(\oX,T\oX)$ gives rise to a bundle map
$$c:\ctX\to T\oX$$
which is an isomorphism outside $x=0$. A curious feature
of the cusp tangent bundle is that it has a canonical \emph{normal} subbundle
to the boundary, rather than conormal as in the standard case. Namely, 
the vector field $x^2\partial_x$ is well-defined (and non-vanishing) 
regardless of changes in the other local coordinates $y_1,\ldots,y_{n-1}$.

\begin{lemma}\label{ldd}
A metric $g_0$ on $X$ extends to a Riemannian metric on the bundle
$\ctX\to \oX$ if and only if it has the form (\ref{cume}). 
\end{lemma}
\dem Clearly, cusp metrics have been defined so that this lemma holds.
\qed

In particular, $g_0$ takes the same form for different local coordinates 
$y_j$, which might not have been \emph{a priori} clear. Notice 
that near $M$ we can write
\begin{equation}\label{fgo}
g_0=a_{00}\frac{dx^2}{x^4}+\frac{dx}{x}\otimes\theta^X+
\theta^X\otimes\frac{dx}{x}+h^X
\end{equation}
where $\theta^X$ and $h^X$ are a smooth $1$-form, respectively a symmetric 
$2$-tensor which restricts to be non-degenerate on the level sets 
$\{x={\rm constant}\}$. There is an ambiguity in the definition
of $\theta^X$, $h^X$ which can be removed by choosing a smooth extension $Z_0$
of the canonical vector field $x^2\partial_x$ to $\oX$ such that
$dx(Z_0)=x^2$, and asking that $\theta^X,h^X$ vanish on $Z_0$.
This vector field corresponds to a product decomposition
of a neighborhood of $M$.

Let now $x'=a x+Bx^2$, with $a\in\cun(M,\rz_+^*)$, $B\in\cun(\oX,\rz)$,
be the Taylor limited development of another boundary-defining function $x'$.
Then $x'$ and $x$ define the same Lie algebra of cusp vector fields
if and only if $a$ is a constant. The set of such functions is called
a cusp structure, which is henceforth fixed. However, it proves useful 
not to fix the function $x$ inconsiderately inside the cusp structure since 
by Lemma \ref{ldd}, the metric $g_0$ is also a cusp metric with respect to $x'$.
The canonical vector fields $x^2\partial_x$ and ${x'}^2\partial_{x'}$ 
are related by $x^2\partial_x=a^{-1}{x'}^2\partial_{x'}$, thus the cusp 
normal bundle is canonically trivialized up to a constant.

Let $\cI$ be the Lie ideal $x\cdot\ccV(\oX)$ inside $\ccV(\oX)$. 
We also denote by $\cI$ the space $x\cun(\oX,\ctsX)$, and more generally
any ideal of the form $x\cun(\oX,E)$ where $E$ is a smooth vector bundle 
over $\oX$.

\begin{lemma}\label{canme}
Let $g_0$ be a cusp metric written in the form (\ref{fgo}).
Let $h:=h^X_{|M}$, $\theta:=\theta^X_{|M}$ and $q:=a_{00}|_M\in\cun(M)$.
Define a $1$-form $\alpha$ and a metric $g^M$ on $M$ by
$$\alpha:=\frac{\theta}{q},\ \ g^M:=\frac{h}{q}-
\alpha\otimes\alpha.$$
Then $q$, $g^M$ and the residue class $\alpha+d\Lambda^0(M)$
are independent, up to a multiplicative constant, 
of the boundary-defining function $x$ inside the 
fixed cusp structure.
\end{lemma}
\dem There is an ambiguity in the choice of $\theta^X$ and $h^X$ but only
up to $\cI^2$, so $\theta$ and $h$ are well-defined.
Clearly, $0<\|x^2-\theta^\#\|_{g_0}^2=q-\|\theta\|_h^2$ so $g^M$ is
non-degenerate. We write
$$g_0=a_{00}\left(\left(\frac{dx}{x^2}
+\frac{\theta^X}{a_{00}}\right)^2+\frac{h^X}{a_{00}}
-\frac{\theta^X}{a_{00}}\otimes\frac{\theta^X}{a_{00}}\right).$$
Let $x'$ be another boundary-defining function in the cusp 
structure of $x$, i.e., $x=ax'+B{x'}^2$ with $a>0$ constant, $B\in\cun(\oX)$. 
Then
$$\frac{dx}{x^2}=\frac{1}{a}\frac{dx'}{{x'}^2}+\frac{dB}{a^2}+\cI$$
so $q'=\frac{q}{a^2}$, $\alpha'=a\alpha+\frac{db}{a}$, where $b=B_{|M}$, and
the lemma follows.
\qed

Melrose \cite{melaps} calls cusp metrics with property (\ref{adias}) exact.  
In light of the above lemma, a cusp metric with $a_{00}{_{|M}}\equiv 1$
and $\theta$ exact can be put in the form (\ref{adias}). By the 
Hodge decomposition, we can always modify $x$ so that $\alpha$
is coexact with respect to $g^M$. 
It makes sense therefore to distinguish closed cusp metrics
as being those with $\theta/q$ closed; in that case, a change of 
boundary-defining function can make $\theta/q$ harmonic with respect to
$g^M$. 

The invariant $\alpha+d\Lambda^0(M)$ has no equivalent in the 
case of finite-length cylinders, see \cite{am}.

\section{$p$- versus $0$-Dirac operators}\label{section3}

Assume that $(X,g_0)$ has a spin structure $\pi:P_\Spin(X)\to P_{SO}(X,g_0)$. 
Let $T_p:(TX,g_p)\to (TX,g_0)$
be the isometry defined by $V\mapsto x^pV$. 
It induces a $SO(n)$-isomorphism
$$T:P_{SO}(X,g_p)\to P_{SO}(X,g_0)$$
between the orthonormal frame bundles of $(X,g_p)$ and $(X,g_0)$.
We define a spin structure on $(X,g_p)$ by 
$T^{-1}\pi:P_\Spin(X)\to P_{SO}(X,g_p)$.
Thus the cusp- and $p$-spinors have the same underlying vector bundle $\Sigma$, 
with the same metrics but different Clifford module structures. The Clifford 
action $c^{g_p}$ with respect to $g_p$ obeys the rule 
$$c^{g_p}(V)\phi=c^{g_0}(T(V))\phi.$$

\begin{prop}\label{isom}
The unbounded operator $D_p$ in $L^2(X,\Sigma,g_p)$ with domain 
$\cun_c(X,\Sigma)$ is isometric to 
$$A_p=x^{-p}\left(D_0-\tfrac{p}{2}xc^0\left(\frac{dx}{x^2}\right)\right)$$
acting in $L^2(X,\Sigma,g_0)$ with domain $\cun_c(X,\Sigma)$.
\end{prop}
\dem 
Note that $dx/x^2$ is a smooth section in $\ctsX$ over $\oX$, thus
$c^0(dx/x^2)$ is actually non-singular down to $x=0$. 

For every conformal transformation $\tg=f^2g$, the
Dirac operators $D_{\tg}$ and $D_g$ are related by \cite[Prop. 1.3]{hitchin}:
\begin{equation}\label{cc}
D_{\tg}=f^{-\tfrac{n+1}{2}}D_g f^{\tfrac{n-1}{2}}.
\end{equation}
Note also the relationship between the two volume densities:
$$dg_p=x^{np} dg_0.$$
Together with the fact that the metric on $\Sigma$ is the same for all $p$,
we see that the map
$$\cun_c(X,\Sigma,g_0)\to \cun_c(X,\Sigma, g_p),\ \ \phi\mapsto 
x^{-\frac{np}{2}}\phi$$
is an isometry with respect to the $L^2$ inner product. Let 
$A_p:=x^{\frac{np}{2}}D_px^{-\frac{np}{2}}$
be the conjugate of $D_p$ under this isometry. 
By using (\ref{cc}) with $f=x^p$, we get 
$A_p=x^{-\frac{p}{2}}D_px^{-\frac{p}{2}}$. The commutation formula
$[D_0,\mu]=c(d\mu)$ for $\mu\in\cun(X)$ ends the proof.
\qed

In particular, the spectra of the closure of $D_p$ and $A_p$ are the same.
It is therefore enough to study $A_p$, which is a cusp operator in the 
sense of the next section.

\section{Review of cusp operators}

\begin{defi}
The algebra $\dc(\oX)$ of cusp differential operators is
the universal 
enveloping algebra of the Lie algebra $\ccV(\oX)$ of cusp vector fields on $\oX$.
\end{defi}

The definition extends easily to 
cusp operators acting on sections of smooth vector bundles over $\oX$.

The spectral properties of cusp differential operators acting as unbounded
operators in $L^2(X, g_0)$
are well-understood, thanks to the calculus of cusp pseudo-differential 
operators of Melrose, a cousin of the celebrated $b$-algebra.
We review below several results about cusp operators, and then use 
them to derive the spectral properties of a class of differential
operators on $L^2(X, g_p)$ which includes the $p$-Dirac operator 
in the spin case. The cusp algebra is described e.g., in 
\cite[Appendix A]{meni96c}. It is a particular case of the fibered cusp algebra
\cite{mame99} when the boundary fibers over a point, and of the cusp algebra
on manifolds with corners \cite{in} when the corners are of 
codimension $1$. 

\subsection{The principal symbol}
Let $A$ be a cusp differential operator. Then the principal symbol of $A$,
which lives \emph{a priori} on the cotangent bundle to $\oX$, lifts in fact
to $\ctsX$ through the (dual) canonical map $T^*\oX\to\ctsX$. We call this 
lifted symbol the \emph{cusp principal symbol}, or simply the principal symbol 
when no confusion can occur. 

\subsection{The normal operator}
This is a "boundary symbol" map, associating 
to any cusp operator $A\in\dc(\oX,\cE,\cF)$ a family of differential
operators on $M$ with one real polynomial parameter $\xi$ as follows:
$$\cN(A)(\xi)=\left(e^{\tfrac{i\xi}{x}}A e^{-\tfrac{i\xi}{x}}\right)|_M$$
where restriction to $M$ is justified by the mapping properties
\begin{eqnarray*}
A:\cun(\oX,\cE)\to\cun(\oX,\cF)&&\\
A:x\cun(\oX,\cE)\to x\cun(\oX,\cF)
\end{eqnarray*}
and by the isomorphism $\cun(M)=\cun(\oX)/x\cun(\oX)$. For example,
$\cN(ix^2\px)(\xi)=-\xi$. The normal operators depends on $x$, nevertheless
if we change $x$ inside its cusp structure, $x'=ax+Bx^2$ with $a>0$ constant and
$B\in\cun(\oX,\rz)$, then the normal operator changes by a conjugation:
$$\cN'(A)(\xi)=e^{-ib\xi/a^2}\cN(A)(\xi/a)e^{ib\xi/a^2}$$
where $b:=B_{|M}$.
In particular, the invertibility of the normal operator for all values
of $\xi\in\rz$ is independent on the particular boundary-defining function $x$.

\subsection{Cusp Sobolev spaces}
For $k\in\nz$, the cusp Sobolev space $H^k_c(\oX,S)$ is defined
as the space of those sections $\phi\in L^2(\oX,S)$ such that $P(\phi)\in
L^2(\oX,S)$ for all $P\in\dc^k(\oX,S)$. It is a Hilbertable space
in which $\cun_c(\oX,S)$ sits as a dense subspace;
moreover, for all $a\in\rz$ an operator $A\in\dc^k(\oX,\cE,\cF)$ has the 
mapping property
$$A:x^aH_c^k(\oX,\cE)\to x^a L^2_c(\oX,\cF)$$
which justifies the definition
$$\dc^{k,a}:=x^{-a}\dc^k,\ \ H_c^{k,a}:=x^a H^{k}.$$

\subsection{The cusp calculus}
As part of his larger program of quantizing singular structures,
Melrose constructed a calculus
$\pc^{z,s}(X)$ of pseudodifferential operators, $z,s\in\cz$, 
in which $\dc$ sits as the differential (i.e., local) operators.
For the construction we refer the reader to \cite{meni96c}.
Operators in $\pc^{z,s}(X)$ extend by duality to distributions and
map $H_c^{k,a}(X)$ into $H_c^{k-\Re(z),a-\Re(s)}(X)$. 
The cusp-principal symbol and the normal operator extend from $\dc(X)$
to this larger calculus. The normal operator takes as values certain families
of pseudodifferential operators on $M$, called 
\emph{suspended operators} \cite{meleta}.

\begin{defi}
An operator $A\in\pc^{z,s}(X)$ is called \emph{fully elliptic} 
if the cusp principal symbol $\sigma_z(x^s A)$ is 
invertible outside the zero section in $\ctsX$, and if moreover 
$\cN(x^sA)(\xi)$ is invertible as
an operator from $H^{\Re(z)}(M)$ to $L^2(M)$ for all values of the 
parameter $\xi$. 
\end{defi}

Fully elliptic operators admit parametrices
with respect to the two symbol structures. By a standard argument we get
\begin{lemma}\label{less}
Let $A\in\pc^{z,s}(X)$ be fully elliptic. Then the domain of the 
closure of $A$ inside $L^2_c(X)$ is $H_c^{\max\{\Re(z),0\},\max\{\Re(s),0\}}(X)$.
If moreover $A$ is symmetric then it is essentially self-adjoint. 
\end{lemma}

The following lemma is the basic tool for the analysis of the spectrum 
of cusp operators. 

\begin{lemma}[Cusp Rellich lemma]
Let $\oX$ be a compact manifold with boundary, $x$ a boundary-defining function
and $g_0$ a cusp metric on $X$.
Then for $p,k\in\rz$ the inclusion $x^pH_c^k(\oX)\hookrightarrow L_c^2(\oX)$ 
is compact if and only if $p>0$ and $k>0$.
\end{lemma}
\dem Easy, using a partition of unity and the classical Rellich lemma 
on a compact exhaustion of a neighborhood of infinity in $X$.
\qed

\begin{cor}[Melrose]
Let $z,s\in\cz$, $a,b\in\rz$ and 
$A\in \pc^{z,s}(\oX,\cE,\cF)$. Then
$A$ is continuous as an operator
$$A:H_c^{a+\Re(z),b+\Re(s)}(\oX,\cE)\to H_c^{a,b}(\oX,\cF),$$
and it is Fredholm if and only if it is fully elliptic.
\end{cor}

\begin{cor} \label{cor14}
The spectrum of a symmetric, fully elliptic operator $A\in\pc^{z,s}(X,S)$ 
is discrete and accumulates towards infinity
if and only if $\Re(z)>0$ and $\Re(s)>0$.
\end{cor}

\section{Complex powers and eigenvalue growth}

Analytic families of cusp operators have been introduced in
\cite{meni96c} to define trace functionals in the spirit of Wodzicki's residue
and to derive an index formula. This idea has been exploited in \cite{in}
for cusp manifolds with corners, and in \cite{phind} for fibered-cusp metrics.
We need here a variant of Melrose' and Nistor's results.

\begin{prop}\label{mert}
Let $p\in\rz_+^*$ and $\cz\supset\Omega \ni z\mapsto A(z)
\in\pc^{z,pz}(X,S)$
an analytic family of cusp operators indexed by a connected open subset 
$\Omega\subset\cz$. Then the trace map
$$\{z\in\Omega;\Re(z)<-n,\Re(pz)<-1\}\ni z\mapsto\Tr(A(z))$$
is holomorphic, and extends meromorphically to $\Omega$ with 
(at most double) possible poles occurring at the superposition
of the sets $(\nz-n)$ and $(\nz-1)/p$.

The first occurring pole is 
\begin{itemize}
\item simple at $z=-n$, if $n>1/p$, with residue 
$$-\tfrac{1}{(2\pi)^n} \int_{\cssX} \tr\sigma_{-n}A(-n);$$
\item double at $z=-n$ if $n=1/p$, with leading coefficient
$$\tfrac{n}{(2\pi)^n} \int_{\cssX_{|M}} \tr\sigma_{-n}(\cN(x^{-1}A(-n)));$$
\item simple at $z=-1/p$, if $n<1/p$, with residue
$$-\tfrac{n}{2\pi} \int_\rz \Tr\cN(x^{-1}A(-1/p))(\xi) d\xi.$$
\end{itemize} 
\end{prop}
\dem (sketch) 
By an analytic family of operators of varying order we mean
that $A$ is holomorphic at $z\in\Omega$ inside the Banach space
of bounded operators
$B(H_c^{k,b},L_c^2)$ for $k>\Re(z)$, $b>p\Re(z)$.

For a trace-class cusp operator $A$ we write its trace as the integral of 
the Schwartz kernel of $A$ on the diagonal in the cusp double-space 
\cite{meni96c}. By the Fourier inversion formula, this becomes the integral
on $\ctsX$ of the pointwise trace of the full symbol of $A$, times the (singular)
canonical volume form given by the cusp symplectic form.
Then we try to extend this integral meromorphically beyond the critical $z$.
This extension follows by a standard elementary argument from the 
polyhomogeneity of the symbol at the boundary of the unit ball bundle, 
see \cite{meni96c} or \cite{in}.
Along the way we collect the information about the leading coefficient in the 
first pole (we could give all asymptotic coefficients in the same way, however 
only the first one matters in the rest of the paper) .
We have omitted the canonical 
volume forms on $\cssX$, resp.\ $\cssX_{|M}$, which are obtained from the symplectic
volume form by contraction with the canonical radial, respectively 
cusp-normal vector field to $M$. 
\qed

\begin{prop}
Let $A\in\pc^{a,b}(X,\cE)$, $a,b>0$, be fully-elliptic, symmetric 
(hence self-adjoint) and positive. Then the complex powers
$A^z$ form an analytic family as in Proposition  \ref{mert}.
\end{prop}
\dem (sketch) This is a particular case of \cite[Thm. 7.1]{phiho}. Here is the 
idea of the proof.
Recall that for classical pseudodifferential operators on a closed manifold
this was shown by Seeley
\cite{seeley}. Guillemin \cite{gui} gave another proof for scalar operators,
which was extended by Bucicovschi to operators on vector bundles, or more 
generally with symbols taking values in a von Neumann algebra. This 
last method extends to 
algebras with "symbols" taking values in algebras like the suspended
algebras, where the existence of complex powers is known.
See also the recent paper \cite{alnv} where complex powers are constructed 
in a more general framework.
\qed

Let $0<\lambda_1\leq\lambda_2\leq\ldots$ be the eigenvalues of $A$. 
The previous two propositions shows that the map
$$z\mapsto \sum_{j=1}^\infty \lambda_j^{-z}=\Tr(A^{-z})$$
is holomorphic for large real parts of $z$, extends meromorphically
to the complex plane with poles on the real axis, and describes the first
pole. To deduce the asymptotic behaviour of $\lambda_j$ we use
a result due to Delange \cite[Theorem III]{delange}:

\begin{lemma}\label{lseq}
Let $k\in\nz^*, a\in\rz^*_+$ and 
$0<\lambda_1\leq\lambda_2\leq\ldots$ be a non-decreasing sequence
such that the series
$$\sum_{j=0}^\infty \lambda_j^{-z}$$
is absolutely convergent for $\Re(z)>a$. Let $\zeta(z)$ denote the 
sum of this series, thus $\zeta(z)$ is holomorphic for $\{\Re(z)>a\}$.
Assume that $\zeta$ extends to a meromorphic function in a neighborhood of 
$\{\Re(z)\geq a\}$, regular except at $z=a$ where it has a pole of order $k$.
Let $C$ be the coefficient of $(z-a)^{-k}$ in the Laurent expansion of
$\zeta$ around $z=a$.
Then the counting function $N(\lambda):=\max\{j\in\nz;\lambda_j<\lambda\}$ 
satisfies
$$N(\lambda)\approx \tfrac{C}{a\cdot (k-1)!}\lambda^a(\log\lambda)^{k-1}.$$
\end{lemma}
\dem Set $\alpha(t):=\sum_{j=1}^\infty H(t-\log \lambda_j)$, where $H$
denotes the Heaviside function. Then
$$f(z):=\int_0^\infty e^{-zt}\alpha(t)dt=\frac{1}{z}\zeta(z)$$
satisfies the hypothesis of \cite[Theorem III]{delange} with leading 
coefficient $C/a$ at $z=a$. Thus 
$$N(e^t)=\alpha(t)\approx \tfrac{C/a}{(k-1)!}e^{at}t^{k-1}.$$
\qed

We note here that 
Karamata's method applied to the zeta function gives a weaker
result than Theorem \ref{th3}, in terms of the asymptotics of the 
spectral function
$$I(n)=\sum_{\lambda_j<n}\frac{1}{\lambda_j}.$$
This is because Karamata's proof does not take into account the behavior of 
the zeta function outside the real line.

We obtain as a corollary a general result about $p$-operators.
\begin{theorem}\label{wg}
Let $X$ be the interior of a compact manifold with boundary $\oX$, with
metric $g_p$ given by (\ref{gp}), $p\in\rz$. Let $\Sigma$ be 
a hermitian vector bundle over $\oX$, $q>0$ and $D\in\pc^{1,q}(X,\Sigma)$ 
a symmetric fully elliptic cusp operator. Then $D$ is essentially 
self-adjoint on $\cun_c(X,\Sigma)$
with discrete spectrum, and the counting function of its eigenvalues satisfies
$$N(\lambda)\approx\left\{\begin{array}{cc}
C\lambda^{n}& \mbox{for $1/n<q$}\\
C\lambda^{n}\log\lambda &\mbox{for $q=1/n$}\\
C\lambda^{1/q}&\mbox{for $0<q<1/n$}.
\end{array}\right.$$
The constant $C$ is explicitly given by Proposition 
\ref{mert} and Lemma \ref{lseq}.
\end{theorem}
\dem $D_q$ is isometric to $A:=x^{np/2}D x^{-np/2}\in
\pc^{1,q}(X,\sigma)$ acting in $L^2(X,\sigma, g_0)$ as in 
Proposition \ref{isom}. Now $\cN(x^{q}A)=\cN(x^qD)$ by the 
multiplicativity of the normal operator. So $A$ is also fully elliptic. The 
theorem follows from Corollary \ref{cor14} and the results of this section.
\qed

Notably, the invariance of the normal operator under conjugation 
by powers of $x$ breaks down for $b$-operators. As a result, the analysis
of conical singularities is more delicate.

\section{The Dirac operator of a cusp metric}\label{docm}

Let us return to the setting of Section \ref{section3}.
Assume that $(X,g_0)$ is spin, and fix a spin 
structure $P_\Spin(X)\to P_{SO}(X)$. Using Lemma \ref{ldd} we extend 
$P_\Spin(X)$ to a spin structure for the Riemannian bundle $(\ctX,g_0)$.
We claim that
the Dirac operator $D_0$ associated to the cusp metric $g_0$ given by
(\ref{cume}) is a cusp differential operator (of order $1$).
For this we first need to describe the smooth structure of the spinor bundle
over $\oX$. The problem is that $g_0$ is singular over the boundary of $\oX$.
Choose a product decomposition 
$$\imath:M\times[0,\epsilon)\to\oX$$
near $M$ such that $x$ becomes the variable in $[0,\epsilon)$, i.e.,
$x(\imath(y,t))=t$. Then $h^X$ (see (\ref{fgo})) defines a Riemannian
metric on the bundle $TM\times[0,\epsilon)$, such that induced map
$\imath:TM\times[0,\epsilon)\to\ctX$ is an isometric injection.

Choose a local orthonormal frame on $M\times\{0\}$, say $Y_1,\ldots,Y_{n-1}$.
Let $Y_1(x),\ldots,Y_{n-1}(x)$ be the orthonormal frame 
obtained by the Gramm-Schmidt procedure from $Y_1,\ldots,Y_{n-1}$, viewed as a
local frame on $M\times\{x\}$.
Let $Y_0$ be inward-pointing unit vector field normal to $M\times\{x\}$ 
for $0<x<\epsilon$. Again by (\ref{cume}), $Y_0$ extends to a smooth 
cusp vector field down to $x=0$. Thus we have found a local orthonormal 
frame in $(\ctX,g_0)$ near $x=0$ which defines a smooth structure on the 
orthonormal frame bundle of $(\ctX,g_0)$, and hence on the spin bundle
$P_{\Spin }(\ctX)$.

An oriented orthonormal frame in $TM\times[0,\epsilon)$ gives rise to a unique
oriented orthonormal frame in $\ctX$ by adding $Y_0$ as the first component.
Thus we get a $SO(n-1)$-map between the principal frame bundles. Define
a spin structure on $M\times\{x\}$, $0\leq x<\epsilon$ as the pull-back of the 
spin bundle of $\oX$ via this map. From the definition there exists
a $\Spin (n-1)$-injection $P_{\Spin }(M\times\{x\})\to P_{\Spin }(\oX)$.

Let $\Sigma_n$ be the spinor representation of $\Spin (n)$. If $n$ is odd,
then $\Sigma_n=\Sigma_{n-1}=\Sigma_{n-1}^+\oplus\Sigma_{n-1}^-$, while 
for $n$ even, $\Sigma_n=\Sigma_{n-1}\oplus\Sigma_{n-1}$. For $n$ odd, 
the Clifford action of 
the extra vector $V_0$ on $\Sigma_{n}$ is given by 
$c^0_{\Sigma_n}=\left[\begin{array}{cc}i&0\\0&-i\end{array}\right]$, and
at the same time, $c^j_{\Sigma_n}=c^0_{\Sigma_n}c^j_{\Sigma_{n-1}}$ 
for $j=1,\ldots,n-1$. For $n$ even, $V_0$ acts as 
$c^0_{\Sigma_n}=\left[\begin{array}{cc}0&-1\\1&0\end{array}\right]$,
while $c^j_{\Sigma_n}=\left[\begin{array}{cc}0&c^j_{\Sigma_{n-1}}\\
c^j_{\Sigma_{n-1}}&0
\end{array}\right]$.

Let $\Sigma M\times\{x\}$, $\Sigma \oX$ be the spinor bundles over 
$M\times \{x\}$, respectively $\oX$. It 
follows from the above discussion that $\Sigma \oX_{|M\times\{x\}}$ 
can be identified with $\Sigma M\times\{x\}$ (for $n$ odd), respectively with
$\rz^2\otimes\Sigma M\times\{x\}$ (for $n$ even) as 
smooth bundles over $\oX$. We wish to compute $D_0$, the Dirac 
operator of the metric $g_0$, in terms of the 
Dirac operators $D^{h^X_x}$ on $M\times\{x\}$. 
First, like for all Dirac operators, $\sigma_2(D_0^2)(\xi)=g_0(\xi,\xi)$.
Lemma \ref{ldd} implies that $D_0^2$ (and hence $D_0$) is cusp-elliptic. 

\begin{lemma}\label{lcco}
The Levi-Civita covariant derivative $\nabla$ on $(X,g_0)$ extends to
a differential operator with smooth coefficients
$$\nabla:\cun(\oX,\ctX)\to \cun(\oX,\ctsX\otimes\ctX).$$
\end{lemma}
\dem Use the Koszul formula for $\nabla$ applied to vectors $Y_j$. Since
$\ctX$ is stable by Lie bracket and $g_0$ is a true metric on
$\ctX$, the assertion follows immediately.
\qed

Let $\tY$ be a local section in the principal 
spin bundle over $\oX$ which lifts the local orthonormal frame 
$(Y_0,\ldots,Y_{n-1})$ constructed previously. Let $\phi$ be a local
section in $\Sigma\oX$. There exists $s:\oX\to \Sigma_n$
smooth so that $\phi=[\tY,s]$, where the square bracket
denotes the equivalence class modulo $\Spin (n)$. In other words, $\phi$
is represented by $s$ in the trivialization given by $\tY$. Recall that
in such a trivialization, 
the Levi-Civita covariant derivative of the spinor $\phi$ takes the form
(at least outside the boundary):
$$\nabla_Y\phi=Y(s)+\tfrac12\sum_{0\leq i<j\leq n-1}
g_0(\nabla_Y Y_i,Y_j)c^ic^j\phi.$$
From Lemma \ref{lcco} it follows that $\nabla$ maps $\cun(\oX,\Sigma\oX)$ to 
$\cun(\oX,\ctsX\otimes \Sigma\oX)$. For $n$ odd, 
\begin{eqnarray}
D_0\phi&=&c^0(D^h\phi+Y_0(s))\nonumber\\&&+\tfrac12\sum_{1\leq i<j\leq n-1}
\langle \nabla_{Y_0} Y_i,Y_j\rangle c^0c^ic^j\phi\nonumber\\
&&+\tfrac12\sum_{1\leq i\leq n-1} \langle \nabla_{Y_i} Y_0,Y_i\rangle c^0\phi
-\tfrac12c^0c(\nabla_{Y_0}Y_0)\phi.\label{dezero}
\end{eqnarray} 
For $n$ even,
\begin{eqnarray*}
D_0^+\phi&=&D^h+
Y_0(s)+\tfrac12\sum_{1\leq i<j\leq n-1}\langle \nabla_{Y_0} Y_i,Y_j\rangle 
c^ic^j\phi\\
&&+\tfrac12\sum_{1\leq i\leq n-1} \langle \nabla_{Y_i} Y_0,Y_i
\rangle \phi -\tfrac12 c(\nabla_{Y_0}Y_0)\phi.
\end{eqnarray*}
and so $D_0\in\dc^1(\oX,\Sigma\oX)$.

We make now the additional hypothesis that
the metric $g_0$ satisfies (\ref{adias}), 
i.e., it is an exact cusp metric \cite{meni96c}.

\begin{lemma} \label{cn}
If $g_0$ is an exact cusp metric then 
\begin{itemize}
\item for $n$ odd, $\cN(D_0)(\xi)=c^0(D^h+i\xi)$;
\item for $n$ even, $\cN(D_0)(\xi)=
\left[\begin{array}{cc}0&D^h-i\xi\\D^h+i\xi&0\end{array}\right]$.
\end{itemize}
Thus $D_0$ is fully elliptic if and only if $D^h$ is invertible.
\end{lemma}
\dem Recall that $\cI$ is the ideal $x\cdot\ccV(\oX)=\ker\cN$ 
inside $\ccV(\oX)$.
From the hypothesis (\ref{adias}), $Y_0=x^2\px+\cI$. It follows that
$\cN(Y_0)=i\xi$ and also $[Y_0,Y_j]\in\cI$ for all $j$.
The Koszul formula shows that $\nabla_{Y_0} Y_i, \nabla_{Y_i} Y_0$
and $\nabla_{Y_0}Y_0$ all belong to $\cI$. Since $\cN(\cI)=0$, 
formula (\ref{dezero}) gives the desired formula.
Then clearly $\cN(D_0)(\xi)$ is invertible for all $\xi\in\rz$ 
if and only if $D^h$ is invertible.
\qed

For example, if $\scal_h(M)\geq 0$ and does not vanish 
identically on any connected component of $M$ then, 
by Lichnerowicz's formula, $D_0$ is fully elliptic.

\section{Proofs of the main results}
In this section (\ref{adias}) is assumed to hold.

{\bf Proof  of Theorem \ref{th1}.}
By Proposition \ref{isom}, we see that for each $p$, $D_p$ is isometric to 
an operator $A_p\in\dc^{1,p}(\oX,\Sigma)$ with the property that 
$A_p-x^{-p}D_0\in\dc^{0,p-1}(\oX,\Sigma)$. 
This means that $A_p$ is fully elliptic if and only $D_0$ is. 
Lemmata \ref{less} and \ref{cn} end the proof.
\qed

{\bf Proof of Theorem \ref{th2}.}
Assume that $D^h$ is invertible, so $D_p\in\pc^{1,p}(\oX,S(\oX))$ 
is fully elliptic by Lemma \ref{cn}. The result follows from 
Corollary \ref{cor14}.
\qed

{\bf Proof of Theorem \ref{th3}.}
We assume that $p>0$ and $D^h$ is invertible. Then $A_p$ is fully elliptic and 
has discrete spectrum, and moreover $P_{\ker A_p}$, the orthogonal
projection onto the kernel of $A_p$ in $L^2(X,\Sigma,g_0)$, belongs to
$\pc^{-\infty,-\infty}(X,\Sigma)$. We apply Proposition \ref{mert} to 
the analytic family
$A(z):=(A_p^2+P_{\ker A_p})^{z/2}\in\pc^{z,pz}(\oX,\Sigma)$ to find that 
the zeta function of $A_p^2$ has either a simple or a double leading pole.
Note that $\sigma_1(A(1))=\|\cdot\|_{g_p}$. From Lemma \ref{cn} 
we get
$$\cN(x^{-1}A(1)^{-1/p})(\xi)=((D^h)^2+\xi^2)^{-\tfrac{1}{2p}},$$
for $n$ odd, respectively
$$\cN(x^{-1}A(1)^{-1/p})(\xi)=\left[\begin{array}{cc}
((D^h)^2+\xi^2)^{-\tfrac{1}{2p}}&0\\0&((D^h)^2+\xi^2)^{-\tfrac{1}{2p}}
\end{array}\right]$$
for $n$ even. Thus 
\begin{eqnarray*}
\tr\sigma_{-n}A(-n)_{|\cssX}&=&\dim\Sigma(n)=2^{[n/2]};\\
\tr\sigma_{-n}(\cN(x^{-1}A(-n)))_{|\cssX_{|M}}&=&\dim\Sigma(n)=2^{[n/2]};\\
\tr\cN(x^{-1}A(-1/p))(\xi)&=&\en\tr((D^h)^2+\xi^2)^{-\tfrac{1}{2p}}.
\end{eqnarray*}
The last identity implies that
\begin{eqnarray*}
\int_\rz \Tr\cN(x^{-1}A(-1/p))(\xi) d\xi&=&\en
\Tr((D^h)^2)^{-\tfrac{1}{2p}+\tfrac{1}{2}}
\int_\rz (1+\xi^2)^{-\tfrac{1}{2p}}d\xi\\
&=&\en\zeta(D^h,\tfrac{1}{p}-1)\tfrac{\sqrt{\pi}\Gamma
\left(\tfrac{1}{2p}-\tfrac{1}{2}\right)}{\Gamma\left(\tfrac{1}{2p}\right)}
\end{eqnarray*}
where the integral was evaluated using \cite[equation (4)]{alim}.
The result follows from these formulae, Proposition \ref{mert} 
and Lemma \ref{lseq}.
\qed

Theorem \ref{th3} applies with minimal modifications to Dirac 
operators twisted by a bundle $W\to\oX$ with 
smooth connection down to $x=0$, provided that the induced twisted 
Dirac operator $D^{M,W}$ over $M$ is invertible. Namely, for
$p\geq 1/n$ we have to multiply by 
$\dim (W)$ the coefficient of $\lambda^n$, resp.\ $\lambda^n\log\lambda$, while
for $0<p<1/n$ we must replace $D^h$ by $D^{M,W}$ inside the zeta function.
If $M$ is connected then the obstruction to the invertibility of
$D^{M,W}$ given by index theory vanishes, as a consequence of 
cobordism invariance. Thus, by Anghel's result \cite{anghel}, 
$D^{M,W}$ is invertible for generic 
connections on $W$, at least if $\dim(M)\leq 4$.

\section{Weyl laws for general cusp metrics}\label{mgen}

Let us consider the case of a general cusp metric. Write $g_0$
in the form (\ref{fgo}) and fix a product decomposition near $M$
so that $\theta^X$ and $h^X$ are uniquely defined.
Let $\theta^\#$ be the dual vector field to $\theta^X$
relative to $h^X$. Set 
$$F:=1/\|x^2\px-\theta^\#\|_{g_0},\ \ f:=F_{|M}\in\cun(M).$$

\begin{prop}\label{norm}
The normal operator of the Dirac operator corresponding to
$g_0$ is given by
\begin{eqnarray}\label{dog}
\cN(D_0)(\xi)&=&c^0(D^h+if\xi-f\nabla_{\theta^\#}\nonumber\\
&&+\tfrac{1}{4} 
fc(d\theta)-\tfrac{1}{4}f\tr( L_{\theta^\#}h(0))-\tfrac{1}{2}
c(df)/f)
\end{eqnarray}
for $n$ odd, respectively
\begin{eqnarray*}
\cN(D_0^+)(\xi)&=&D^h+if\xi-f\nabla_{\theta^\#}\\
&&+\tfrac{1}{4} 
fc(d\theta)-\tfrac{1}{4}f\tr( L_{\theta^\#}h(0))-\tfrac{1}{2}
c(df)/f
\end{eqnarray*}
for $n$ even.
\end{prop}
\dem The idea is to expand formula (\ref{dezero}) for $D_0$. Note that
$Y_0=F(x^2\px-\theta^\#)$, and therefore 
\begin{eqnarray}
[Y_i,Y_0]&=&\frac{Y_i(F)}{F}Y_0+F(-[Y_i,\theta^\#]-x^2\px(Y_i))\nonumber\\
&\equiv&\frac{Y_i(f)}{f}Y_0-f[Y_i,\theta^\#]+\cI^2\label{lb}
\end{eqnarray}
(recall that $\cI$ is the ideal $x\ccV(\oX)$ annulated by $\cN$). Thus
\begin{eqnarray*}
\cN(D_0)(\xi)&=&c^0\cdot(D^h+if\xi-f\nabla_{\theta^\#})\\
&&+\tfrac12\sum_{1\leq i<j\leq n-1}
\langle \nabla_{fx^2\px} Y_i,Y_j\rangle_{|x=0} c^0c^ic^j\\
&&+\tfrac12\sum_{1\leq i\leq n-1} \langle \nabla_{Y_i} Y_0,Y_i\rangle_{|x=0} c^0\\
&&-\tfrac12\sum_{j=1}^{n-1} \langle \nabla_{Y_0}Y_0,Y_j\rangle_{|x=0} c^0c^j.
\end{eqnarray*}
We already see here the first three terms from (\ref{dog}). 
Koszul's formula gives 
$$2\langle \nabla_{x^2\px} Y_i,Y_j\rangle=
Y_i(\theta_j)-Y_j(\theta_i)-\theta([Y_i,Y_j])=d\theta(Y_i,Y_j)$$
so by summing over $i<j$ we get the fourth term. Since $Y_i$ has constant length
$1$, it follows by (\ref{lb}) that 
$$\langle \nabla_{Y_i} Y_0,Y_i\rangle=
\langle[Y_i,Y_0], Y_i\rangle =-f\langle [Y_i,\theta^\#],Y_i\rangle
+\cI^2.$$
Use 
$$\sum_{i=1}^{n-1}\langle [Y_i,\theta^\#],Y_i\rangle=
\tfrac12\tr(L_{\theta^\#}h(0))=\frac{L_{\theta^\#}d h(0)}{d h(0)}$$
to get the fifth term.
As for the last term, again by (\ref{lb}) we have
$$\langle \nabla_{Y_0}Y_0,Y_j\rangle=-\langle Y_0,\nabla_{Y_0}Y_j\rangle=
\langle Y_0,[Y_j,Y_0]\rangle\equiv Y_j(f)/f+\cI^2.$$
This settles the case $n$ odd. The case $n$ even is done in the same way.
\qed

We wish now to give conditions on $g_0$ so that $\cN(D_0)(\xi)$ is invertible
for all $\xi\in\rz$. If $\theta=0$ this condition was seen to be 
equivalent to $D^{M,h(0)}$ being invertible. For the general case, 
we rewrite (\ref{dog}) as
\begin{eqnarray*}
\lefteqn{\cN(D_0)(\xi)=}\\&&=c^0f\left(f^{-\frac12}D^hf^{-\frac12}+i\xi
-\left(\nabla_{\theta^\#}+\tfrac{1}{2}\frac{L_{\theta^\#}d h(0)}{d h(0)}\right)
+\tfrac{1}{4}c(d\theta)\right).
\end{eqnarray*}
It is not clear if a reasonable 
condition on $h,f$ and $\theta$ exists so that
the above operator is invertible for all $\xi\in\rz$. Recall that the metric
$h$ is not canonically determined by $g_0$, since it depends on the choice of
$x$ inside the fixed cusp structure. Thus, we rewrite $\cN(D_0)$ in terms of
the canonical metric $g^M$ on $M$, see Lemma \ref{canme}, although
we must still fix $x$ in order to define $\cN$.

\begin{prop}\label{normd}
Let $g_0$ be a cusp metric on $X$ and $\alpha,g^M,q$ defined in Lemma
\ref{canme}.
Then for $n$ odd,
$$\cN(D_0)(\xi)=c^0q^{-\frac{n+1}{2}}(D^{g^M}+i\xi
(1-c(\alpha))+\tfrac14 c(d\alpha))q^{-\frac{n-1}{2}}$$
while for $n$ even,
$$\cN(D_0^+)(\xi)=q^{-\frac{n+1}{2}}(D^{g^M}+i\xi
(1-c(\alpha))+\tfrac14 c(d\alpha))q^{-\frac{n-1}{2}}.$$
\end{prop}
\dem
Recall that in the expression (\ref{fgo}) for $g_0$, the tensors
$\theta^X,h^X$ are well-defined only
up to $\cI^2$. Fix a product decomposition of $\oX$ near $M$, thus 
removing the ambiguity. 

First, assume that ${a_{00}}=1$ near $M$, so $g_0=(dx/x^2+\theta^X)^2+h^X-
\theta^X\otimes\theta^X$. Let $\tY_0:=x^2\px$.
Define an isometric embedding
$$(TM\times[0,\infty),h^X)\to(\ctX,g_0),\ \ Y\mapsto \tY:=Y-\theta(Y)\tY_0.$$ 
This allows us to compare spinors on $X$ and on $M$. Choose a local orthonormal
frame $Y_1,\ldots,Y_{n-1}$ on $M$, transport it to $X$ using the 
product decomposition and re-orthonormalize it using Gramm-Schmidt. Note that
the frame $\{Y_j\}$ is different from the frame with the same name 
from Section \ref{docm}, 
since they are orthonormal with respect to different metrics.
Then $\tY_0,\ldots,\tY_{n-1}$ is a (smooth) orthonormal frame on $\ctX$.
We use $i,j,k$ to denote a subscript in $\{1,\ldots,n-1\}$.
Notice that
\begin{eqnarray*} [\tY_i,\tY_j]&=&\widetilde{[Y_i,Y_j]}-
d\theta(Y_i,Y_j)+\cI^2\\ 
{[}\tY_0,\tY_i]&\in&\cI^2.
\end{eqnarray*}
So 
\begin{eqnarray*}
(\nabla_{\tY_i}\tY_j,\tY_k)&=&(\nabla_{Y_i}Y_j,Y_k)\\
\nabla_{\tY_0}\tY_0&\in&\cI^2\\
(\nabla_{\tY_i}\tY_0,\tY_i)&\in&\cI^2\\
(\nabla_{\tY_0}\tY_j,\tY_k)&=&\tfrac12 d\theta(\tY_i,\tY_j)+\cI^2.
\end{eqnarray*}
Using (\ref{dezero}), under the assumption that
$a_{00}=1$ we get immediately in the case $n$ is odd
$$\cN(D_0)(\xi)=c^0(D^{g^M}+i\xi(1-c(\theta))+\tfrac14 c(d \theta))).$$

Remove now the assumption that $a_{00}=1$ in a neighborhood of $M$.
Using the formula for the conformal change for 
Dirac operators (\ref{cc}), we get the desired
expression. The case $n$ even is entirely similar.
\qed

As a corollary, we give a condition for the full ellipticity
of $D_p$ in some cases where $g_0$ is not exact.

\begin{cor}\label{cor22}
Let $g_0$ be a cusp metric with $\alpha$ closed (by Lemma \ref{canme}, this
condition is independent on the boundary-defining function).
Then $D_p$ is fully elliptic if and only if $D^{g^M}$ is invertible. 
\end{cor}
\dem By Proposition \ref{normd}, we must check the invertibility of the
family of operators
$$P(\xi):=D^{g^M}-i\xi c(\alpha)+i\xi$$
for all $\xi\in\rz$. But $P(\xi)$ is invertible if and only if 
$P^*(\xi)P(\xi)=(D^{g^M}-i\xi c(\alpha))^2+\xi^2$ is invertible.
This holds automatically for $\xi\neq 0$, and for $\xi=0$ it 
is equivalent to $D^{g^M}$ being invertible.
\qed

Product cylindrical metrics have been studied by several authors starting with
\cite{aps1}. Such metrics can be considered either 
as $b$- or as cusp metrics; they can actually be 
treated by elementary methods using separation of variables.
Melrose \cite{melaps} has successfully studied exact $b$-metrics, which are
in some sense only asymptotically cylindrical. 
Corollary \ref{cor22} allows even more geometric structure to be embedded in
the metric. If $H^1(M,\rz)\neq 0$, it simply says that our results on Weyl laws
hold for some metrics (i.e., closed cusp metrics) which are not exact. 

Let us finally state our most general result about Dirac eigenvalues of cusp 
metrics.

\begin{theorem}
Assume that the normal operator of $D_0$, computed in two different
ways in Propositions \ref{norm}, \ref{normd}, is invertible. Then
for $p>0$, $D_p$ is essentially self-adjoint with pure point spectrum
accumulating towards infinity. The rate of growth is given by Theorem \ref{th3}
if $p\geq 1/n$, and by Theorem \ref{wg} if $0<p<1/n$. 
\end{theorem}
\qed

This theorem has the drawback that the invertibility of $\cN(D_0)$ 
must be assumed. The only case where we have found a reasonable
condition for invertibility is for closed cusp metrics (Corollary \ref{cor22}). 
We are therefore let to the following

{\bf Problem.} Let $(M,g^M)$ be a connected compact Riemannian spin manifold
and $\theta$ a $1$-form on $M$. Find conditions
on $\theta, g^M$ so that the family of operators 
$$D^{g^M}+i\xi(1-c(\theta))+\tfrac14 c(d \theta))$$
is invertible for all $\xi\in\rz$.

Another variant of the problem would be to assume additionally that 
$M$ vanishes in the Spin bordism ring. Mere invertibility of $D^{g^M}$
is not enough in general, as shown by the following example (courtesy
of Andrei Moroianu, see also \cite{am}):

\begin{example}
Let $M=S^2$ with an arbitrary metric $h$ 
and $D$ the associated Dirac operator. Then $D^h$ is invertible, since
the Dirac operator of the standard metric is invertible, and any two metrics
on $S^2$ are conformally isometric. Nevertheless, there exist a metric $h$
and $\theta\in\Lambda^1(S^2)$ so that $D^h+c(\theta)$ is not invertible.
For this, take $h$ to be the metric induced from an immersion of $S^2$ 
in $\rz^3$ with total mean curvature $0$. First, there exists such an 
immersion, since we can deform $S^2$ by stretching cylinders of positive 
or negative mean curvature, thus increasing or decreasing the total 
mean curvature at will. 
Secondly, take $\tilde{\phi}$ to be a constant spinor on
$\rz^3$, and pull it back to a spinor $\phi=(\phi^+,\phi^-)$ on $M$. Then
by \cite[Proposition 2]{friedrich},
$$D^h\phi^+=-iH\phi,\ \ \ \ D^h\phi^-=i\phi^+$$
where $H$ is the mean curvature function of $(M,h)$. 
The form $Hdh$ has volume $0$, thus it is exact, i.e., $Hdh=d\theta$ for some
$\theta\in\Lambda^1(M,\rz)$. Recall moreover that $c(dh)$ acts by $\pm i$
on $\Sigma^\pm(M)$.
Thus $\phi$ is a solution of the equation
$$(D^h+c(d\theta))\phi=0.$$
\end{example}

\bibliographystyle{plain}

\end{document}